\documentclass[12pt]{article}[smallcondensed]
\usepackage{pgf,tikz}
\usetikzlibrary{arrows}
\usepackage[USenglish]{babel}
\usepackage{amsmath,amsfonts,amssymb,amsthm}
\usepackage[utf8]{inputenc}
\usepackage[all]{xy}
\usepackage[T1]{fontenc}
\usepackage{xfrac}
\usepackage{authblk}
\usepackage{mathtools}
\usepackage{pifont}
\usepackage{marvosym}

\usepackage{geometry}
\usepackage{tikz}
\usetikzlibrary{matrix}
\usepackage{tikz-cd}
\newtheorem{teo}{Theorem}[section]
\newtheorem{defi}[teo]{Definition}
\newtheorem{prop}[teo]{Proposition}

\newtheorem{coro}[teo]{Corollary}

\newtheorem*{ejem}{Example}
\theoremstyle{definition}
\newtheorem*{dem}{Proof}

\newtheorem*{cond}{Condition}

\providecommand{\keywords}[1]
{
  \small	
  \textbf{\textit{Keywords---}} #1
}

\title{\textbf{A category theory approach using preradicals to model information flows in networks}}

\author[1,3]{\small Sebastian Pardo G}
\author[1,2,3 \Letter]{Gabriel A. Silva}

\affil[1]{Department of Bioengineering\\University of California, San Diego}
\affil[2]{Department of Neurosciences\\University of California, San Diego}
\affil[3]{Center for Engineered Natural Intelligence\\University of California, San Diego}

\date{}

\begin{document}
\maketitle

\renewenvironment{abstract}
{\begin{quote}
\noindent \rule{\linewidth}{.5pt}\par{\bfseries \abstractname.}}
{\medskip\noindent \rule{\linewidth}{.5pt}
\end{quote}
}
\begin{abstract}
    Category theory has been recently used as a tool for constructing and modeling an information flow framework. Here, we show that the flow of information can be described using preradicals. We prove that preradicals generalize the notion of persistence in spaces where the underlying structure forms a directed acyclic graph. We show that a particular $\alpha$ preradical describes the persistence of a commutative $G$-module associated with a directed acyclic graph. Furthermore, given how preradicals are defined, they are able to preserve the modeled system's underlying structure. This allows us to generalize the notions of standard persistence, zigzag persistence, and multidirectional persistence. \\
\end{abstract}

Corresponding author: GS, gsilva@ucsd.edu\\

\keywords{Information flow, preradicals, quiver representation, persistence module, persistence.}

\section{Introduction}

This work is the first of two related papers that aim to explore a formal framework for the conceptual axiom that, in a general sense, information flow through a network will always follow a path of 'least resistance'. We propose that the resultant effect of such a least resistance information flow is ultimately responsible for the dynamics on the network, independent of the physical processes that carry the information in any specific network. While distinct physical processes are responsible for driving the flux of information and determining what the unit or carrier of information is in a specific network, a least resistive information flow path is a universal and generalized phenomenon. The critical intuition here is that while descriptive models or algorithmic rules set up the conditions for the dynamics of a network, it is the flow of signals and information through the network - bounded by physical constraints - that dictate how computational events are prioritized and executed. This naturally follows a path of least resistance. Yet, to the best of our knowledge, a theoretical framework that captures the universality of this effect and allows computing such information flow least resistive pathways in real networks does not exist. Here, we attempt to provide the theoretical foundation for such a framework.

For example, in a literal way, water flowing through a set of interconnected pipes will flow down a path of least resistance as dictated by the pressure fronts associated with the flow. The physics responsible for this includes the incompressibility of water and the compliance and diameter of the pipes in the network. The resultant effect is a 'natural' flow determined and constrained by the physics that takes the route of least effort. 

Similarly, weight changes in artificial neural networks are shaped by constraints. Their  physics is in the connectivity structure of the network, error or loss functions, and the algorithms (rules) that affect weight changes, e.g., backpropagation. The consequence of these constraints and rules are decentralized but integrated weight changes that allow the network to 'learn'. The weight changes themselves, however, follow a flow of least effort or resistance that is a natural consequence of the physical make up and rules of the network. 

Another example are biological neural networks composed of interconnected neurons. In fact, this was our main practical motivation for this work. There are myriad of quantitative models that describe different aspects of neural dynamics and the structural rules that specify how brain networks are connected. Yet, describing how internal computations and information flows are prioritized is not currently possible. We do not fully understand the foundational principles that govern how networks of neurons successfully operate autonomously. What the drivers are that shape the prioritization, order, and decisions of the internal computations and operations after a neuron and network of neurons receive inputs? All this is distinct from specifying an algorithm a network needs to follow or a model of the neuron. Mechanistically, in the case of biological neural networks, an example of a fundamental physical processes that impose such constraints are structure (connectivity and geometry)-function constraints, and energy considerations \cite{Gabe}. 

Likewise, in human social networks, we tend to make friends with individuals we relate to and like after introductions are made. Network connectivity is one physical constraint that sets the condition for the creation of pathways that reflect a least effort or resistance, i.e., individuals we like. In many systems, including artificial neural networks and  biological brain networks, these processes are responsible for passing information from one scale of organization or hierarchy to another. 

Here, we use category theory to describe how information flows through a system that can be modeled using objects and morphisms in the category $K$-Mod, where $K$ is a field. We show that the flow of information is preserved within the modeled system's underlying structure. We achieve this by using a theory that involves persistent modules, which derive from persistence homology. Persistent homology considers a family of topological spaces and inclusions from one space into another to find topological features common to subsets within these spaces. The persistent homology is defined by homology groups that allow classifying invariant topological features. We can think, for example, of the standard persistent homology where spaces are linearly ordered with one space included in the next. Similarly, we can think of zigzag persistence \cite{Zig Zag}, where the spaces are linearly ordered, but the inclusions can occur in either direction. Another example is multidimensional persistence \cite{Multi Pers}, which operates in multiple dimensions on a grid with inclusion maps parallel to the coordinate axis. 

In \cite{Pers} Chambers and Letsher extended these notions of persistence by considering that the underlying structure forms a directed acyclic graph (DAG), naming it DAG persistence. The authors used an algebraic representation called $G$-\textit{module} endowed with an additional commutative condition that arises from following distinct paths in a DAG between the same pair of vertices. The construction of this representation is based on algebraic structures which can be considered objects in the category of $K$-modules for $K$ a field. In this paper, we will describe persistence from a category theory perspective, using a tool known as preradicals. 

Preradicals produce compatible substructures within each object preserved by morphisms, allowing for an efficient determination of algebraic invariants. Using preradicals we were able to show how information flows through the $G$-module representation. We will generalize the notion of persistence by showing that the persistence of a $G$-module, as defined in \cite{Pers}, can be obtained employing a particular type of preradical in the category $K$-Mod. This generalization enables us to use preradicals to examine how information is transmitted between the objects that comprise the $G$-module representation so that the underlying structure is preserved.  

The paper is organized as follows: in Section \ref{Pers module} we present preliminary definitions for directed acyclic graphs and quivers. We then give the basic construction of a quiver representation and the definition of persistence in a $G$-module. In Section \ref{Prerradicals} we introduce the definition of preradicals along with some notation, terminology, and basic facts. We discuss the role of preradicals as a tool to describe how information is transmitted within the category. We conclude this section with the definition of the $\alpha$ and $\omega$ preradicals. Our main contributions are introduced in Section \ref{ Persistence}, where we show that the persistence of a $G$-module whose underlying structure is a single-source single-sink graph, is obtained by an $\alpha$  preradical. We then generalize this proposition to a $G$-module whose underlying structure defines a graph with $n$ sources and $m$ sinks. Section \ref{Interpretations} discuss preradicals as a general way to describe the flow of information in a quiver representation modeled in the category $R$-Mod for $R$ a commutative ring with unity. Section \ref{Conclusion} provides some concluding arguments.  In Appendix \ref{Apendix} we provide a brief description of the construction of homology groups as well as the persistence homology group.

\section{Preliminaries} \label{Pers module}

We begin by giving preliminary definitions of relevance to the rest of our paper. For a complete introduction to homology groups, we direct the reader to some of the standard references \cite{Bredon}, \cite{Hatcher}, \cite{Munkers} and \cite{Rotman}. Most models of persistence use a collection of spaces and inclusions from one space into another to find common topological features among these spaces. In \cite{Pers} the authors show a generalization that considers inclusions over a set of spaces whose underlying structure forms a directed graph, with the constraint that the graph must be a directed acyclic graph (DAG): acyclic and not contain repeated edges. We call this structure, comprised of the collection of spaces with inclusions from one space into another, a graph filtration. Formally, 
\begin{defi}
For a simple directed acyclic graph $G=(V,E)$, a graph filtration $\chi_{G}$ of a topological space $X$ is a pair $(\{X_{v}\}_{v\in V}, \{f_{e}\}_{e\in E})$ such that
\begin{itemize}
    \item[(1)] $X_{v}\subset X$ for all $v\in V$;
    \item[(2)] If $e=(v,u)\in E$ then $f_{e}:X_{v}\longrightarrow X_{u}$ is a continuous embedding (or inclusion) of $X_{v}$ into $X_{u}$.
\end{itemize}
\end{defi}

A \textit{quiver} is a directed graph where loops and multiple directed edges between the same vertices are allowed. Formally, a quiver is a quadruple $Q=(V,E,s,t)$ where $V$ represents the set of \textit{vertices}, $E$ the set of \textit{edges} and $s,t:E\longrightarrow V$ are two maps, assigning the \textit{starting vertex} and the \textit{ending vertex} for each edge. In this case, for an edge $e$ with $s(e)=u$ and $t(e)=v$ we write $e:u\longrightarrow v$. A quiver $Q$ is finite if both sets $V$ and $E$ are finite.  In particular, every DAG is a quiver. 

Typically quivers admit a \textit{representation} over a field $K$ which assigns to each vertex $v$ a vector space $W_{v}$ and to each arrow $e:u\longrightarrow v$ a $K$-linear morphism $f_{e}:W_{u}\longrightarrow W_{v}$. More formally,
\begin{defi}
Given a quiver $Q=(V,E,s,t)$, a representation of $Q$ over a field $K$ is a pair of families
\begin{center} 
$\Big(\{W_{v}\}_{v\in V}, \{f_{e}\}_{e \in E} \Big)$  
\end{center}
where for each edge $e:u\longrightarrow v$, $f_{e}:W_{u}\longrightarrow W_{v}$ is a $K$-linear morphism.
\end{defi}

Since we are interested in quivers with relations, precisely commutative conditions, we will add the following:

\begin{cond}
All diagrams defined by two directed paths in the quiver $Q$ commute: in other words, for any path $\gamma=e_{1},\cdots, e_{n}$ in $Q$, one can extend this to a $K$-linear morphism $f_{\gamma}=f_{e_{n}}\circ \cdots \circ f_{e_{1}}$ in the quiver representation. Then, commutativity means that for any two different directed paths $\gamma$ and $\gamma'$ in $Q$ connecting vertices $u$ and $v$, we have that $f_{\gamma}=f_{\gamma'}$.   
\end{cond}

\begin{defi}\cite[Definition 2.2]{Pers}
For a DAG $G=(V,E)$, a \textbf{commutative $G$-module} is the pair $\Big(\{W_{v}\}_{v\in V}, \{f_{e}\}_{e \in E} \Big)$ where for each vertex $v$, $W_{v}$ is a vector space and for any edge $e=(u,v)$, $f_{e}:W_{v}\longrightarrow W_{u}$ is a $K$-linear morphism with the condition that the resulting diagrams are commutative. 
\end{defi}

To put it another way, since every DAG is a quiver, then a commutative $G$-module for a DAG $G=(V,E)$ is the same as a quiver representation with commutative conditions. The next definition provides the framework to generalize the notion of standard persistence and multidirectional persistence, which will be shown in Section \ref{ Persistence}.

\begin{defi}\cite[Definition 2.3]{Pers} For a directed acyclic graph $G=(V,E)$ a $k$-dimensional persistence module for a graph filtration $\chi_{
G}$, is the commutative $G$-module $(\{W_{v}\}_{v\in V}, \{f_{e}\}_{ e\in E})$ where 
\begin{itemize}
    \item[(i)] $W_{v}=H_{k}(X_{v})$\footnote{Here $H_{k}(X_{v})$ denotes the $k$-homology group of $X_{v}$. See Appendix.} for all $v\in V$;
    \item[(ii)] For every edge $e=(u,v)\in E$, $f_{e}:H_{k}(X_{u})\longrightarrow H_{k}(X_{v})$ is the map induced by the the inclusion $X_{u}\hookrightarrow X_{v}$. 
\end{itemize}
\end{defi}

Suppose now we are given a directed acyclic graph $G=(V,E)$ and its commutative $G$-module representation, denoted as $M=(\{W_{v}\}_{v\in V}, \{f_{e}\}_{ e\in E})$. A cone in the category $K$-Mod for the $G$-module $M$ is a pair $(L,\{\eta_{u}\}_{u\in V})$ where $L$ is a vector space and $\eta_{u}$ is a $K$-linear morphism for each $u\in V$, such that for any edge $e=(u,v)$, we have that the diagram
\[
\begin{tikzcd}
[column sep=small]
& L \ar[dl, "\eta_{u}"'] \arrow[dr, "\eta_{v}"] & \\
W_{u}\arrow{rr}{f_{e}} &  & W_{v},
\end{tikzcd}
\]
commutes in $K$-Mod, that is, $f_{e}\circ \eta_{u}=\eta_{v}$. Thus, a limit for the $G$-module $M$ is a cone $(L,\{\eta_{u}\}_{u\in V})$ with the universal property: for any other cone $(L',\{\eta'_{u}\}_{u\in V})$ there exist a unique morphism $\eta:L'\longrightarrow L$ such that $\eta_{u}\circ \eta = \eta'_{u}$ for every $u\in V$:
\[
\begin{tikzcd}
[column sep=small]
 & L' \ar[dl, "\eta"'] \arrow[dr, "\eta'_{u}"] & \\
L \arrow{rr}{\eta_{u}} &  & W_{u}.
\end{tikzcd}
\]

By duality, we define a cocone for the $G$-module $M$ as a pair $(L, \{\eta_{u}\}_{u\in V})$ where $L$ is a vector space and $\eta_{u}$ is a $K$-linear morphisms for each $u\in V$, such that for any edge $e=(u,v)$, the following diagram commutes in $K$-Mod

\[
\begin{tikzcd}
[column sep=small]
W_{u} \arrow{rr}{f_{e}} \ar[dr, "\eta_{u}"'] &    & W_{v} \ar[dl, "\eta_{v}"]\\
 & L. & 
\end{tikzcd}
\]

A colimit for the $G$-module $M$ is a cocone $(L, \{\eta_{u}\}_{u\in V})$ with the universal property: for any other cocone $(L', \{\eta'_{u}\}_{u\in V})$ there exist a unique morphism $\eta:L\longrightarrow L'$ such that $ \forall u\in V$, the diagram 
\[
\begin{tikzcd}
[column sep=small]
 & W_{u} \ar[dl, "\eta_{u}"'] \arrow[dr, "\eta'_{u}"] & \\
L \arrow{rr}{\eta} &  & L'
\end{tikzcd}
\]
is commutative in $K$-Mod.

The category $K$-Mod is known to be complete and cocomplete, that is, where all small limits and colimits exist. Since small categories index the commutative diagrams into consideration, then the limit and colimit always exist for a commutative $G$-module. Therefore, for any commutative $G$-module $M$, we have an induced morphism
\begin{center}
    $\varphi_{M}:\textsl{lim}(M)\longrightarrow \textsl{colim}(M)$. 
\end{center}
This is the precursor to defining the persistence of a $G$-module.

\begin{defi}(Definition 2.5 in \cite{Pers}.)
If $M$ is a commutative $G$-module then the persistence of $M$, denoted by $P(M)$, is the image $\varphi_{M}(\textsl{lim}(M))$.
\end{defi}

The persistence $P(M)$ of a $G$-module $M$ tells us how much information of $\textsl{lim}(M)$ persists through the $G$-module structure but not about the evolvement of the information that persists from each component of the $G$-module. It is in this context that preradicals can provide information about the evolution of the structure $\varphi_{M}(\textsl{lim}(M))$, as well as the substructures from each component throughout the $G$-module representation. On the one hand, by applying a preradical $\sigma$ to $\textsl{lim}(M)\overset{\varphi_{M}}{\longrightarrow}\textsl{colim}(M)$ we get how much information defined by $\sigma$ is preserved through the $G$-module representation. On the other hand, by applying a preradical $\sigma$ to each component of the $G$-module representation, we obtain a sketch of how the information defined by $\sigma$ changes through the entire $G$-module structure.

\section{Preradicals}\label{Prerradicals}

This section introduces preradicals and a sketch of how information is transmitted in a diagram whose components are comprised of objects and morphisms in the category of $R$-Mod. For a complete introduction to preradicals and its properties in $R$-Mod see \cite{Bican}, \cite{Pre 1} and \cite{Partitions}. As introduced in \cite{Isom}, preradicals can also be considered as compatible choice assignments, since these assign to each object in the category a subobject, in such a way that the restriction and corestriction of any morphism to the corresponding subobjects, results in a morphism within the category. 

\begin{ejem}
Consider the category $\mathbb{Z}$-Mod of all $\mathbb{Z}$-modules. This category is isomorphic to the category $Ab$ of all abelian groups. Now, given $M\in \mathbb{Z}$-Mod, one can define   
\begin{center}
$\sigma(M)=\{x\in M \mid 2x=0\}$.
\end{center}

Notice first that $\sigma(M)$ is a submodule of $M$. Also, for any morphism $f:M\longrightarrow M'$ in $\mathbb{Z}$-Mod and any $y\in M$, we have that $f(2y)=2f(y)$. Thus, if $x\in \sigma(M)$ we have that $2f(x)=f(2x)=0$, which implies that $f(x)\in \sigma(M')$.
Hence, $f(\sigma(M))\subseteq \sigma(M')$ which in turns implies that
\begin{equation*}
\xymatrix{
M\ar[r]^{f} & M' \\
\sigma(M) \ar[r]^{f_{|}} \ar[u]^{\iota} & \sigma(M')\ar[u]_{\iota}
}
\end{equation*}
is a commutative diagram. 
\end{ejem}

The above example shows an endofunctor \footnote{An endofunctor is a functor whose domain is equal to its codomain.} on $\mathbb{Z}$-Mod which acts as a subfunctor of the identity functor on $\mathbb{Z}$-Mod. This is precisely the definition of a preradical:

\begin{defi}
Let $\mathcal{C}$ be a category. A preradical $\sigma$ on the category $\mathcal{C}$ is a functor  that assigns to each object $C\in \mathcal{C}$, a subobject $\sigma(C)$ such that for each morphism $f:C\longrightarrow C'$ in $\mathcal{C}$, we have the commutative diagram 
\begin{equation*}
\xymatrix{
C\ar[r]^{f} & C' \\
\sigma(C) \ar[r]^{\sigma(f)} \ar[u]^{\iota} & \sigma(C')\ar[u]_{\iota}
.}
\end{equation*}
Here, $\sigma(f)$ is the restriction and corestriction of $f$ to $\sigma(C)$ and $\sigma(C')$ respectively, that is, $\sigma(f):=f\mid_{\sigma(C)}:\sigma(C)\longrightarrow \sigma(C')$. Also, $\iota$ represents the inclusion map.
\end{defi}

One can define operations between preradicals for certain categories. For example, in the category of left $R$-modules (where $R$ is a commutative unitary ring) as well as in the category of complete modular lattices $\mathcal{L}_{\mathcal{M}}$, one can define four principal operations (see \cite{Pre 1} and \cite{Isom} respectively). For the purpose of this work, we will consider the category $R$-Mod since we will use quiver representations as a model to study the flow of information defined by a directed graph. 
\\ We denote by $R$-pr the collection of all preradicals on $R$-Mod. There is a natural partial ordering in $R$-pr given by $\sigma \leq \tau$ if and only if $\sigma(M)\leq \tau(M)$ for all $M\in R$-Mod. The four classical operations in $R$-pr are defined as follows: if $M\in R$-Mod and $\sigma,\tau\in R$-pr, then

\begin{itemize}
\item[i)] $(\sigma \wedge \tau)(M)=\sigma(M)\cap \sigma(M)$;
\item[ii)] $(\sigma \vee \tau)(M)= \sigma(M) + \tau(M)$;
\item[iii)] $(\sigma \cdot \tau)(M)=\sigma (\tau(M))$;
\item[(iv)] $(\sigma:\tau)(M)$ is such that 
\begin{center}
$(\sigma:\tau)(M)/\sigma(M)=\tau(M/\sigma(M))$.
\end{center} 
\end{itemize}
The above operations are called the meet, the join, the product and the coproduct respectively. Now, as the intersection and the sum of an arbitrary family of $R$-modules defines an $R$-module, then both the join and the meet operations can be defined for any family of preradicals, thus making $R$-pr a big lattice \footnote{A big lattice is a class (not necessarily a set) having joins and meets for arbitrary families (indexed by a class) of elements.}: for any family $\{\tau_{\alpha}\}_{\alpha \in I}$ of preradicals and $M$ an $R$-module, 
\begin{itemize}
\item[i)] $\big( \bigvee_{\alpha\in I} \ \tau_{\alpha} \big)(M) = \Sigma_{\alpha\in I} \ \tau_{\alpha}(M)$,

\item[ii)] $\big( \bigwedge_{\alpha\in I} \ \tau_{\alpha} \big)(M) = \cap_{\alpha\in I} \ \tau_{\alpha}(M)$.
\end{itemize}

From now on, we will assume that the base ring $K$ is a field, although the results in this section hold for any commutative ring with unit. A \textit{diagram of shape} $\mathcal{J}$ in the category $K$-Mod is a functor $F:\mathcal{J}\longrightarrow K\textit{-}Mod$. Commonly, we use \textit{small} categories to define diagrams in a category, since these have only a set's worth of arrows (and thus, of objects as well since every object in the category defines an identity morphisms).
For instance, if we have a linear ordered shape diagram (that is, where $J$ is a poset) 
\begin{equation*}
\xymatrix{
\cdots \ar[r] & M_{i}\ar[r]^{f} & M_{i+1}\ar[r]^{g} & \dots \ar[r] & M_{j} \ar[r]^{h} & M_{j+1} \ar[r] & \cdots  }
\end{equation*}
we can apply any preradical $\sigma$, obtaining the commutative diagrams 
\begin{equation*}
\xymatrix{
 \cdots \ar[r] & M_{i}\ar[r]^{f} & M_{i+1}\ar[r]^{g} & \dots \ar[r] & M_{j} \ar[r]^{h} & M_{j+1} \ar[r] & \cdots   \\
\cdots \ar[u]^{\iota} \ar[r] & \sigma(M_{i}) \ar[r]^{f_{|}} \ar[u]^{\iota} & \sigma(M_{i+1}) \ar[r]^{g_{|}} \ar[u]^{\iota} & \dots \ar[r] \ar[u]_{\iota} & \sigma(M_{j})\ar[u]_{\iota} \ar[r]^{h_{|}} & \sigma(M_{j+1}) \ar[u]_{\iota} \ar[r] & \cdots \ar[u]^{\iota} 
.}
\end{equation*}
Since all diagrams are commutative, the collection $\{\sigma(M_{i})\}_{i\in J}$ sketches how the information, regarding $\sigma$, flows through the diagram of shape $J$. Also, if we consider that the collection of morphisms with their possible compositions provide the information contained in the diagram of shape $J$, then preradicals--as compatible choice assignments--show how these morphisms transmit the information within the diagram. Moreover, the four operations defined on $K$-pr enable us to study the evolvement of the information through the system. For example, the meet of a family $\{\sigma_{i}\}_{i\in I}$ of preradicals can be used to get the common information to all $\sigma_{i}$ within each component $M$ of the diagram, whereas the join of the family $\{\sigma_{i}\}_{i\in I}$ of preradicals will give us the full information generated by all $\sigma_{i}$, within $M$. Likewise, from the product $\sigma \cdot \tau$ we obtain the information defined by $\sigma$ after the information by $\tau$ is processed; while the coproduct $(\sigma : \tau)$ gives the information of $\tau$ conditioned to contain the information defined by $\sigma$.

Finally, we define in $K$-Mod the preradicals $\alpha$ and $\omega$ as follows: given a $K$-module $M$ and a submodule $N$ of $M$, the preradical $\alpha_{N}^{M}$ evaluated in a $K$-module $W$ is 
\begin{center}
$\alpha_{N}^{M}(W)= \sum \{f(N)\mid f\in Hom(M,W)\}$,
\end{center}
where $Hom(M,W)=\{f:M\longrightarrow W)\}$. Likewise, given $N$ a submodule of $M$, the preradical $\omega_{N}^{M}$ evaluated in a $K$-module $W$ is
\begin{center}
$\omega_{N}^{M}(W)= \cap \{f^{-1}(N)\mid f\in Hom(W,M)\}$,
\end{center}
where $Hom(W,M)=\{f:W\longrightarrow M)\}$ and $f^{-1}(N)$ denotes the inverse image of $N$ under the morphism $f$. One can easily prove that every preradical $\sigma$ in $K$-Mod can be written as 
\begin{center}
    $\sigma = \bigvee \{\alpha_{\sigma(M)}^{M} | M\in K\textit{-}Mod\}=\bigwedge \{\omega_{\sigma(M)}^{M} | M\in K\textit{-}Mod \}.$
\end{center}

We end this section with the next result which will be used in the discussions of Section \ref{Interpretations}:
\begin{prop}\cite{Bican}[Proposition I.1.2] \label{preradicales suma directa}
Let $\sigma$ be a preradical and $\{M_{i}\}_{i\in I}$ be a family of $R$-modules. Then $\sigma(\oplus_{i\in I} \ M_{i})=\oplus_{i\in I}\ \sigma(M_{i})$.
\end{prop} 

\section{Persistence Through Preradicals} \label{ Persistence}

In this section we discuss our main technical contributions and results. Let us first notice that given any quiver $Q$, one can generate the \textit{path category} of the quiver $Q$, which we denote as $\Lambda(Q)$. The collection of objects of the category $\Lambda(Q)$ coincides with the collection of vertices of $Q$, whereas the collection of morphisms in $\Lambda(Q)$ is comprised by all finite paths between two vertices of the quiver. Here, the composition operation is given by concatenation of paths. Thus, we can think of a commutative $G$-module representation as a diagram of shape $\Lambda(Q)$ in $K$-Mod whose commutative condition holds. This diagram is given by the functor $F:\Lambda(Q)\longrightarrow K\textit{-}Mod$ induced by the quiver representation, which assigns to each object in $\Lambda(Q)$ a $K$-module and to each path in $\Lambda(Q)$ a $K$-linear morphism. This fact will allow us to apply the $\alpha$ preradicals on $K$-Mod to express the flow of information within the $G$-module representation associated to a directed acyclic graph $G=(V,E)$. To this end, the morphisms to be considered in each $\alpha_{N}^{M}$ will be restricted to the ones that appear in the $G$-module representation. This will be achieved by composing the functor $F$ mentioned above, followed by an $\alpha$ preradical into consideration. 

Suppose now we have a single\textit{-}source single\textit{-}sink graph $G$ and its commutative $G$-module representation $M=(\{W_{v}\}_{v\in V}, \{f_{e}\}_{ e\in E})$
\begin{equation*}
\xymatrix{
W_{s}\ar[r]^{f_{s,t}} & W_{i}\ar[r]^{f_{i,k}} & \dots \ar[r] & W_{j} \ar[r]^{f_{j,t}} & W_{t}
,}
\end{equation*}
where $W_{s}$ is the $K$-module representing the source vertex while $W_{t}$ is the $K$-module representing the sink vertex of the graph $G$. As its shown in \cite[Lemma 3.3]{Pers}, the limit and the colimit of the commutative $G$-module $M$ is given by $W_{s}$ and $W_{t}$ respectively. Thus, the persistence of $M$ is 
\begin{center}
    $P(M)=im(W_{s}\overset{\varphi}{\longrightarrow} W_{t})=\varphi(W_{s})$,
\end{center}
where $\varphi$ is the induced $K$-linear morphism from $lim(M)$ to $colim(M)$. As we next show, the persistence $P(M)$ in \cite[Lemma 3.3]{Pers} can also be obtained using an $\alpha$ preradical: 

\begin{prop}\label{Prop 1}
Let $M=(\{W_{v}\}_{v\in V}, \{f_{e}\}_{ e\in E})$ be a commutative $G$-module such that the underlying directed acyclic graph is a single\textit{-}source single\textit{-}sink graph, with $W_{s}$ and $W_{t}$ representing the source and the sink vertices in $M$, respectively. Then the persistence shown in \cite[Lemma 3.3]{Pers}, is the same as the preradical $\alpha_{W_{s}}^{W_{s}}$ evaluated in $W_{t}$. In this case, the morphisms $g:W_{s}\longrightarrow W_{t}$ that occur in the description of the preradical $\alpha_{W_{s}}^{W_{s}}$ are considered within the $G$-module representation.
\end{prop}

\begin{dem}
Let $M=(\{W_{v}\}_{v\in V}, \{f_{e}\}_{ e\in E})$ be a commutative $G$-module such that the underlying directed acyclic graph is a single\textit{-}source single\textit{-}sink graph. Let us denote by $W_{s}$ and $W_{t}$ the representation of the source and the sink vertices in $M$, respectively. By definition, the preradical $\alpha_{W_{s}}^{W_{s}}$ evaluated at $W_{t}$ is 
\begin{center}
    $\alpha_{W_{s}}^{W_{s}}(W_{t})=\sum \Big \{g(W_{s})\mid g\in Hom(W_{s},W_{t})\Big \}$.
\end{center}
Here, the set $Hom(W_{s},W_{t})$ is taken as all $K$-linear morphisms in $K$-Mod with domain $W_{s}$ and codomain $W_{t}$. However, as noted at the beginning of this section, we can restrict the morphisms in $Hom(W_{s},W_{t})$ to the ones present in the commutative $G$-module representation, still obtaining commutative diagrams
\begin{equation*}
\xymatrix{
W_{s}\ar[r]^{f_{s,i}} & W_{i}\ar[r]^{f_{i,k}} & \dots \ar[r] & W_{j} \ar[r]^{f_{j,t}} & W_{t}  \\
\alpha_{W_{s}}^{W_{s}}(W_{s}) \ar[r]^{f_{s,i}|} \ar[u]^{\iota} & \alpha_{W_{s}}^{W_{s}}(W_{i}) \ar[r]^{f_{i,k}|} \ar[u]^{\iota} & \dots \ar[r] & \alpha_{W_{s}}^{W_{s}}(W_{j})\ar[u]_{\iota} \ar[r]^{f_{j,t}|} & \alpha_{W_{s}}^{W_{s}}(W_{t}) \ar[u]_{\iota}
.}
\end{equation*}
As the underlying graph is a single-source single-sink graph, the limit and the colimit of its commutative $G$-module representation $M$ are the objects $W_{s}$ and $W_{t}$, respectively. Thus, since the persistence of $M$ is given by the image of the induced morphism $\varphi:W_{s}\longrightarrow W_{t}$, the commutative condition of the $G$-module representation implies that  $g(W_{s})=\varphi(W_{s})$ for every morphism $g\in Hom(W_{s},W_{t})$ in the $G$-module representation $M$. Thus,
\begin{align*}
    \alpha_{W_{s}}^{W_{s}}(W_{t})=\sum \Big \{g(W_{s})\mid g\in M \mbox{ and } g\in Hom(W_{s},W_{t})\Big \} \\
    =\sum \varphi(W_{s}) =\varphi(W_{s}) =im(\varphi)=P(M).
\end{align*}
\qed
\end{dem}

From the preceding proposition, we see that for any other preradical $\sigma$ in $K$-Mod, when applied to each component of the diagram that comprise the $G$-module representation, we also obtain commutative diagrams
\begin{equation*}
\xymatrix{
W_{s}\ar[r]^{f} & W_{i}\ar[r]^{g} & \dots \ar[r] & W_{j} \ar[r]^{h} & W_{t}  \\
\sigma(W_{s}) \ar[r]^{f_{|}} \ar[u]^{\iota} & \sigma(W_{i}) \ar[r]^{g_{|}} \ar[u]^{\iota} & \dots \ar[r] & \sigma(W_{j})\ar[u]_{\iota} \ar[r]^{h_{|}} & \sigma(W_{t}) \ar[u]_{\iota}
.}
\end{equation*}
Since each diagram commutes and the $G$-module $M$ possesses the commutative condition, it follows that the information defined by the preradical $\sigma$ that persists through the $G$-module representation is $\varphi(\sigma(W_{s}))\subseteq im(\varphi)\cap \sigma(W_{t})$. Furthermore, we can apply any of the four operations $\sigma \wedge \tau$, $\sigma \vee \tau$, $\sigma \cdot \tau$ and $(\sigma:\tau)$ between preradicals to study the persistence of information within the $G$-module representation. Here, we can see how the commutative diagrams induced by preradicals helps us represent the persistence and thus the flow of information: they preserve the underlying structure of the diagram. 

As shown in \cite[Proposition 3.4]{Pers}, DAG persistence generalizes the Standard persistence, the Zigzag persistence, and the Multidimensional persistence. In their proof, the authors rely on \cite[Lemma 3.3]{Pers} to prove that the Standard and Multidimensional persistence homology groups can be defined in terms of the image of a map. In contrast, for zigzag persistence, the definition of the zigzag persistence module is identical to the definition of DAG persistence module. In the following result, we show the corresponding proposition using preradicals assuming that the graph filtration also satisfies the commutative diagram condition. Again, the morphisms considered in the $\alpha$ preradical are those which appear in the commutative $G$-module representation.

\begin{prop}
Suppose $\chi_{G}$ is a graph filtration of $X$ satisfying the commutative diagram condition and let $\Big( \{H_{k}(X_{v})\}_{v\in V}, \{f_{e}\}_{e\in E}\Big)$ its commutative $G$-module representation. Then:
\begin{itemize}
    \item[(1)] (Standard Persistence) \\ If $G$ is the graph corresponding to the filtration $X_{0}\rightarrow X_{1} \rightarrow \cdots \rightarrow X_{n}$ and $I_{i,p}$ is the subgraph consisting of vertices $\{X_{i}, \cdots , X_{p}\}$ then $H_{k}^{I_{i,p}}(\chi_{G})= \alpha_{H_{k}(X_{i})}^{H_{k}(X_{i})}(H_{k}(X_{i+p}))$.
    
    \item[(2)] (Multidimensional Persistence) \\ Let $\chi =\{X_{v}\}_{v\in \{0,\dots,m\}^{d}}$ be a multifiltration with underlying graph $G$. If $G_{u,v}$ is the subgraph with vertices $\{w \in G | u \leq w \leq v \}$ then the rank invariant $\rho_{X,k}(u,v)=dim \ \alpha_{H_{k}(X_{u})}^{H_{k}(X_{u})}(H_{k}(X_{v}))$
\end{itemize}
\end{prop}

\begin{dem}
(1) From the fact that $H_{k}^{I_{i,p}}(\chi_{G})=im(H_{k}(X_{i})\longrightarrow H_{k}(X_{i+p}))$ and that the graph filtration satisfies the commutative diagram condition, it follows that 
\begin{center}
  $H_{k}^{I_{i,p}}(\chi_{G})=im(H_{k}(X_{i})\longrightarrow H_{k}(X_{i+p}))=\alpha_{H_{k}(X_{i})}^{H_{k}(X_{i})}(H_{k}(X_{i+p}))$.  
\end{center}
(2) The rank invariant $\rho_{X,k}(u,v)$ is defined as the dimension of the image of the induced map $H_{k}(X_{u})\longrightarrow H_{k}(X_{v})$. Since the graph filtration satisfies the commutative diagram condition, then any two morphisms in the $G$-module representation associated with distinct paths from $u$ to $v$ will have the same image. Hence, we have that 
\begin{center}
   $H_{k}^{G_{u,v}}=im(H_{k}(X_{u})\longrightarrow H_{k}(X_{v}))=\alpha_{H_{k}(X_{u})}^{H_{k}(X_{u})}(H_{k}(X_{v}))$.
\end{center}
Thus, the rank invariant is the dimension of the subspace defined by the preradical $$\alpha_{H_{k}(X_{u})}^{H_{k}(X_{u})}(H_{k}(X_{v})).$$
\qed
\end{dem}

We will now generalize Proposition \ref{Prop 1} by assuming that the underlying DAG have $n$ sources $s_{1},\dots, s_{n}$ and $m$ sinks $t_{1},\dots, t_{m}$. Before we do this, we shall make the following observation: suppose that $M=(\{W_{v}\}_{v\in V}, \{f_{e}\}_{ e\in E})$ is a commutative $G$-module representation for a given DAG $G=(V,E)$ having $n$ sources $s_{1},\dots,s_{n}$ and $m$ sinks $t_{1},\dots, t_{m}$:

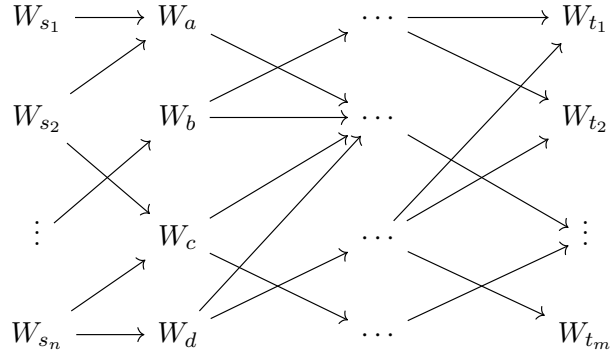
\begin{center}
\begin{tikzcd}
W_{s_{1}} \arrow[r]             & W_{a} \arrow[rrd]              &  & \cdots \arrow[rr] \arrow[rrd]               &  & W_{t_{1}} \\
W_{s_{2}} \arrow[ru] \arrow[rd] & W_{b} \arrow[rr] \arrow[rru]   &  & \cdots \arrow[rrd]                          &  & W_{t_{2}} \\
\vdots \arrow[ru]               & W_{c} \arrow[rru] \arrow[rrd]  &  & \cdots \arrow[rruu] \arrow[rrd] \arrow[rru] &  & \vdots    \\
W_{s_{n}} \arrow[r] \arrow[ru]  & W_{d} \arrow[rruu] \arrow[rru] &  & \cdots \arrow[rru]                          &  & W_{t_{m}}
\end{tikzcd}
\end{center}

Then, for any cone $\Big(C,\{\eta_{v}\}_{v\in V}\Big)$ of $M$ we have that all morphisms $\eta_{j}:C\longrightarrow W_{j}$ with $j\notin \{s_{1},\cdots, s_{n}\}$, are factorized by some $\eta_{s_{i}}:C\longrightarrow W_{s_{i}}$ with $i\in \{1,\dots, n\}$. Indeed, when considering the diagram induced by the cone 
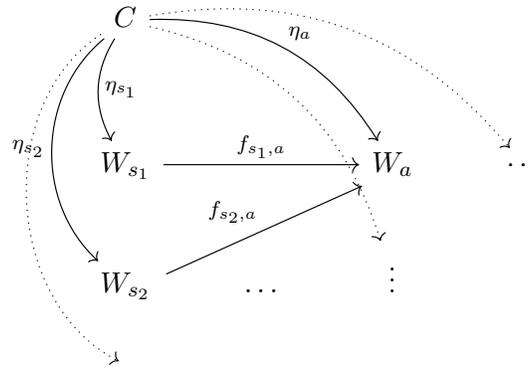
\begin{center}
\begin{tikzcd}
C \arrow[rrdd, "\eta_{a}", bend left] \arrow[dd, "\eta_{s_{1}}", bend right] \arrow[ddd, "\eta_{s_{2}}"', bend right=49] \arrow[rrrdd, dotted, bend left] \arrow[rrddd, dotted, bend left] \arrow[dddd, dotted, bend right=60] &       &        &        \\
                                                                                                                                                                                                                               &       &        &        \\
W_{s_{1}} \arrow[rr, "{f_{s_{1},a}}"]                                                                                                                                                                                              &       & W_{a}  & \cdots \\
W_{s_{2}} \arrow[rru, "{f_{s_{2},a}}", shift right]                                                                                                                                                                                & \dots & \vdots &        \\
{}                                                                                                                                                                                                                             &       &        &       
\end{tikzcd}
\end{center}
the commutative property defined by $\Big(C,\{\eta_{v}\}_{v\in V}\Big)$ implies that
\begin{equation*}
 f_{s_{1},a}\circ \eta_{s_{1}}=\eta_{a}=f_{s_{2},a}\circ \eta_{s_{2}}. 
\end{equation*}
In other words, for each $l\neq \{s_{1},\cdots, s_{n}\}$ the morphism $\eta_{l}$ is factorized by some $\eta_{s_{i}}$, with $i=1,\cdots, n$. Furthermore, the above argument shows that any two or more factorizations of $\eta_{l}$ are equal. Thus, when obtaining the persistence of the commutative $G$-module $M$, it suffices to take the object $lim(M)$ along with the respective morphisms $\eta_{s_{i}}:lim(M)\longrightarrow W_{s_{i}}$ for $i=1,\cdots,n$ in the cone $\Big(lim(M),\{\eta_{v}\}_{v\in V}\Big)$ of $M$. A similar argument shows that we only need to consider the object $colim(M)$ together with the respective morphisms $\iota_{t_{i}}:W_{t_{i}}\longrightarrow colim(M)$, for $i=1,\cdots, m$ in the cocone $\Big(colim(M),\{\iota_{v}\}_{v\in V}\Big)$ when computing the persistence of the commutative $G$-module $M$. 
\\Lastly, to define the $\alpha$ preradical, which will be used to get the persistence, we need to define a submodule of $lim(M)$ that contains in some way the information that flows from the $n$ sources. This will be given by the submodule generated by the inverse images $\eta^{-1}_{s_{1}}(W_{s_{1}}),\cdots, \eta^{-1}_{s_{n}}(W_{s_{n}})$. In fact, one can verify that all these inverse images coincide: let $\Big(C,\{\eta_{v}\}_{v\in V}\Big)$ be any cone of the commutative $G$-module $M$ and let us take the product $\Big(\Pi_{i=1}^{n} \ W_{s_{i}}, \{\rho_{i}\}_{i=1}^{i=n}\Big)$ of the source objects $W_{s_{1}},\cdots W_{s_{n}}$, where each $\rho_{i}$ represents the projection on $W_{s_{i}}$. Then, the universal property of the product induces the following diagram
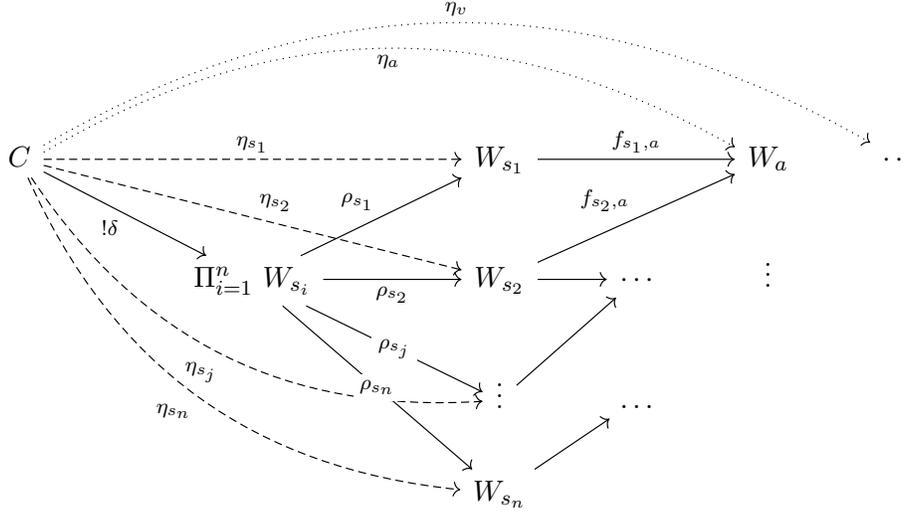
\begin{center}
\begin{tikzcd}
C \arrow[rrd, "!\delta"'] \arrow[rrrr, "\eta_{s_{1}}", dashed] \arrow[rrrrd, "\eta_{s_{2}}", dashed] \arrow[rrrrddd, "\eta_{s_{n}}"', dashed, bend right] \arrow[rrrrrr, "\eta_{a}"', dotted, bend left, shift right] \arrow[rrrrrrr, "\eta_{v}", dotted, bend left] \arrow[rrrrdd, "\eta_{s_{j}}"', dashed, bend right] &  &                                                                                                                                                                    &  & W_{s_{1}} \arrow[rr, "{f_{s_{1},a}}"]             &        & W_{a}  & \cdots \\
                                                                                                                                                                                                                                                                                                                         &  & \Pi_{i=1}^{n} \ W_{s_{i}} \arrow[rru, "\rho_{s_{1}}"] \arrow[rrdd, "\rho_{s_{n}}" description] \arrow[rr, "\rho_{s_{2}}"'] \arrow[rrd, "\rho_{s_{j}}" description] &  & W_{s_{2}}  \arrow[rru, "{f_{s_{2},a}}"] \arrow[r] & \cdots & \vdots &        \\
                                                                                                                                                                                                                                                                                                                         &  &                                                                                                                                                                    &  & \vdots \arrow[ru]                                 & \dots  &        &        \\
                                                                                                                                                                                                                                                                                                                         &  &                                                                                                                                                                    &  & W_{s_{n}} \arrow[ru]                              &        &        &       
\end{tikzcd}
\end{center}
which implies that $\eta_{s_{i}}=\rho_{s_{i}} \circ \delta$ for all $i\in \{1,\dots, n\}$. Since each $\rho_{s_{i}}$ is a surjective morphism, we have that 
\begin{center}
    $\eta_{s_{i}}^{-1}(W_{s_{i}})=\big(\rho_{s_{i}} \circ \delta\big)^{-1}(W_{s_{i}}) =\delta^{-1}\big(\rho_{s_{i}}^{-1}(W_{s_{i}})\big)=\delta^{-1}\Big(\Pi_{i=1}^{n} \ W_{s_{i}}\Big)$,
\end{center}
for all $i\in \{1,\cdots, n\}$. Therefore, the submodule of $lim(M)$ used in the definition of the $\alpha$ preradical to describe persistence will be $\underbar{M}:=\delta^{-1}\Big(\Pi_{i=1}^{n} \ W_{s_{i}}\Big)$. 

\begin{prop} \label{Prop lim colim}
Let $M=(\{W_{v}\}_{v\in V}, \{f_{e}\}_{ e\in E})$ be a commutative $G$-module such that the underlying directed acyclic graph have $n$ sources and $m$ sinks. If $W_{s_{1}},\cdots W_{s_{n}}$ represent the source nodes in the commutative $G$-module representation, then the persistence of the $G$-module is obtained by evaluating the preradical
\begin{center}
     $\alpha_{\underbar{M}}^{lim(M)}(colim(M))=\sum \Big \{f(\underbar{M}) \ |\ f: lim(M) \longrightarrow colim(M) \Big \}$
\end{center}
where the morphisms $f: lim(M) \longrightarrow colim(M)$, considered in the description of the preradical $\alpha_{\underbar{M}}^{lim(M)}$, are taken within the extended $G$-module representation by the limit and the colimit of $M$.
\end{prop}

\begin{dem}
Let $M$ be a commutative $G$-module such that the underlying directed acyclic graph have $n$ sources and $m$ sinks and let $W_{s_{1}},\cdots W_{s_{n}}$ represent the $n$ source nodes in the $G$-module representation. From the observations made previously, we obtain an extension of the commutative $G$-module representation as a single\textit{-}source single\textit{-}sink diagram:
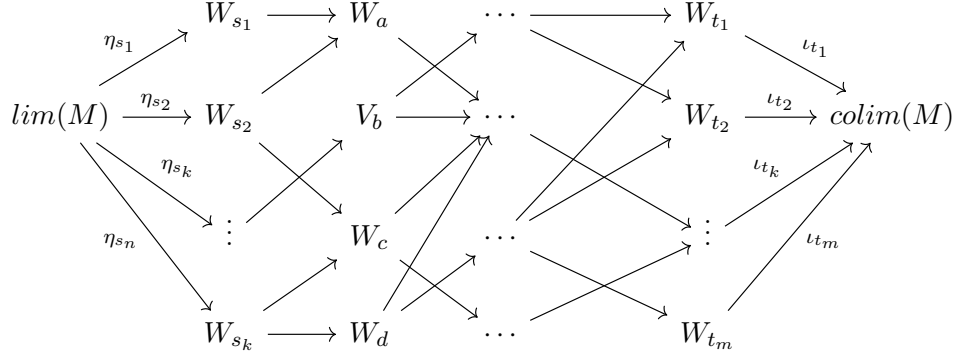
\begin{center}
\begin{tikzcd}
                                                                                                                                & W_{s_{1}} \arrow[r]             & W_{a} \arrow[rd]             & \cdots \arrow[rr] \arrow[rrd]               &  & W_{t_{1}} \arrow[rd, "\iota_{t_{1}}"]   &          \\
lim(M) \arrow[ru, "\eta_{s_{1}}", shift left] \arrow[r, "\eta_{s_{2}}"] \arrow[rd, "\eta_{s_{k}}"] \arrow[rdd, "\eta_{s_{n}}"'] & W_{s_{2}} \arrow[ru] \arrow[rd] & V_{b} \arrow[r] \arrow[ru]   & \cdots \arrow[rrd]                          &  & W_{t_{2}} \arrow[r, "\iota_{t_{2}}"]    & colim(M) \\
                                                                                                                                & \vdots \arrow[ru]               & W_{c} \arrow[ru] \arrow[rd]  & \cdots \arrow[rruu] \arrow[rrd] \arrow[rru] &  & \vdots \arrow[ru, "\iota_{t_{k}}"]      &          \\
                                                                                                                                & W_{s_{k}} \arrow[r] \arrow[ru]  & W_{d} \arrow[ruu] \arrow[ru] & \cdots \arrow[rru]                          &  & W_{t_{m}} \arrow[ruu, "\iota_{t_{m}}"'] &         
\end{tikzcd}
\end{center}

Observe that since every morphism $f:lim(M)\longrightarrow colim(M)$ is factorized by some $\eta_{s_{i}}$, then when $f$ is restricted to $\underbar{M}$ we assure that the information flowing from $lim(M)$ is passing through some $W_{s_{i}}$, for $i=1,\dots, n$. Therefore, it follows that the persistence of the commutative $G$-module $M$ is given by
\begin{center}
     $\alpha_{\underbar{M}}^{lim(M)}(colim(M))=\sum \Big \{f(\underbar{M}) \ |\ f: lim(M) \longrightarrow colim(M) \Big \}$
\end{center}
where all morphisms $f: lim(M) \longrightarrow colim(M)$ are taken within the extended $G$-module representation.
\qed
\end{dem}

We can also use $\alpha$ preradicals to describe the persistence of information defined by a proper subset of sinks. In this case, for the sources $W_{s_{i_{1}}},\cdots W_{s_{i_{l}}}$ with $0\leq l<n$, we must consider the submodule $\underbar{M}$ of $lim(M)$ together with all morphisms $f:lim(M)\longrightarrow colim(M)$ which are factorized by some $\eta_{s_{i_{j}}}$, with $j\in \{1,\dots, l\}$.  
\begin{coro}
Let $M$ be a commutative $G$-module such that the underlying directed acyclic graph have $n$ sources and $m$ sinks. The persistence of the sources $W_{s_{i_{1}}},\cdots W_{s_{i_{l}}}$ for $0\leq l<n$ is obtained by evaluating the preradical
\begin{center}
     $\alpha_{\underbar{M}}^{lim(M)}(colim(M))=\sum \Big \{f(\underbar{M}) \ |\ f: lim(M) \longrightarrow colim(M)\Big \}$
\end{center}
where all morphisms $f:lim(M)\longrightarrow colim(M)$ are factorized by some $\eta_{s_{i_{j}}}$ with $j\in \{1,\dots, l\}$.
\end{coro}

In a sense, the $\alpha_{N}^{M}$ preradicals fulfill the role of describing the propagation of information forward-wise. In other words, when evaluating the preradical $\alpha_{N}^{M}$ in $K$ we obtain the \textit{direct} information in $K$ coming from the submodule $N$ within the domain $M$ (the morphisms go from $M \longrightarrow K$). If in the quiver representation we find morphisms between components going in either direction, we can make use of the $\omega$ preradicals to express \textit{indirect} information. Here, the morphisms considered in the definition of the $\omega$ preradicals go from $K\longrightarrow M$, thus obtaining the information in $K$ induced by the inverse images of the submodule $N$ which lies within the codomain $M$. This indirect information also persists through the quiver representation and is obtained similarly, since any preradical is a compatible choice assignment that induces commutative diagrams: given a $K$-module $W_{j}$ and a submodule $L_{j}$, the following diagram commutes
\begin{equation*}
\xymatrix{
lim(M)\ar[r]^{\varphi} & colim(M) \\
\omega_{L_{j}}^{W_{j}}(lim(M)) \ar[r]^{\varphi_{|}} \ar[u]^{\iota} & \omega_{L_{j}}^{W_{j}}(colim(M))\ar[u]_{\iota}
.}
\end{equation*}

\section{Preradicals and their Interpretation }\label{Interpretations}

In this last section, we sketch how preradicals describe the flow of information in a system modeled with objects and morphisms in the category $R$-Mod. Motivated by principles from Shannon's version of information theory, where messages are created by choosing letters or words from a dictionary, preradicals are also functions of choice whose choices are compatible with the structure of the category. This property gives preradicals a high order approximation to the information contained in a  system modeled in $R$-Mod. We will approach this from a local perspective, although the techniques employed here also describe the flow of information through any set of objects and morphisms in $R$-Mod. The latter follows from the fact that preradicals are functors that preserve the underlying structure of any diagram. 

Given a system, for example, a neural network, we can decompose it into a collection of objects called components related to each other by an input-output relationship. These input-output relationships connect the components and can be represented by an oriented arrow $A\longrightarrow B$ if an output of $A$ is an input of $B$, giving the system a directed graph representation. However, such a representation cannot describe many situations of the actual input-output existing relations between the components of the system. For instance, a component $A$ may contribute to multiple distinct inputs, although the directed graph representation will only show one arrow connecting $A$ to $B$. 
To overcome this, we say that two components $A$ and $B$ are connected if there exists a function $f:A\longrightarrow B$; this allows us to consider the set of all different input-output relations between these two components. Observe that each set of functions describing the input-output relations is determined by the components $A$ and $B$. As we have mentioned before, any directed graph can be associated with a category $\Lambda(Q)$ whose objects represent the vertices of the graph and whose morphisms correspond to finite paths between vertices in the graph. This association of a directed graph to its category of paths leads us to define a diagram of shape $\Lambda(Q)$ in $R$-Mod via a functor $F:\Lambda(Q)\longrightarrow R$-Mod.

Let us now consider a diagram of shape $\Lambda(Q)$ in $R$-Mod and let $B$ be a component in this diagram that receives two inputs from components $A_{1}$ and $A_{2}$. Then, we have two morphisms $f_{1}:A_{1}\longrightarrow B$ and $f_{2}:A_{2}\longrightarrow B$ which represents an input-output relationship between the respective source and target components. Both functions define a transformation from the direct sum (or direct product since it has a finite number of components) $f:A_{1}\oplus A_{2}\longrightarrow B$ such that $f\circ \iota_{i}=f_{i}$, that is,
\begin{center}
$\big( A_{i}\overset{\iota_{i}}{\hookrightarrow} A_{1}\oplus A_{2} \overset{f}{\longrightarrow} B\big) =A_{i}\overset{f_{i}}{\longrightarrow} B$,
\end{center}
where $\iota_{i}:A_{i}\longrightarrow A_{1}\oplus A_{2}$ denotes the inclusion map.
\\ If $\sigma$ is a preradical in $R$-Mod, then by Proposition \ref{preradicales suma directa} we have $\sigma(A_{1}\oplus A_{2})=\sigma(A_{1})\oplus \sigma(A_{2})$, and thus the diagram
\begin{equation*}
\xymatrix{
A_{1}\oplus A_{2}\ar[r]^{\ \ \ \ f} & B \\
\sigma(A_{1})\oplus \sigma(A_{2}) \ar[r]^{\ \ \ \ \ f|} \ar[u]^{\iota} & \sigma(B)\ar[u]_{\iota}
}
\end{equation*}
commutes. With this in mind, we obtain the information defined by the compatible choice assignment $\sigma$ relative to the transformation $f$ within the system. This information turns out to be 
\begin{center}
$f(\sigma(A_{1})) + f(\sigma(A_{2}))\subseteq \sigma(B)\subseteq B$.
\end{center}

We can generalize the previous argument to a number of $n$ inputs from components $A_{1},\cdots, A_{n}$. In this case, we get a set of $n$ morphisms 
\begin{center}
$A_{1}\overset{f_{1}}{\longrightarrow} B$,\\
$ \vdots$ \\ 
$A_{n}\overset{f_{n}}{\longrightarrow} B$,
\end{center}
which induces a morphism $f:\oplus_{i=1}^{n} A_{i}\longrightarrow B$ such that $f\circ \iota_{i}=f_{i}$, that is,
\begin{center}
$\big( A_{i}\overset{\iota_{i}}{\hookrightarrow} \oplus_{i=1}^{n} A_{i} \overset{f}{\longrightarrow} B\big) =A_{i}\overset{f_{i}}{\longrightarrow} B$, for each $i\in \{1,\cdots, n\}$.
\end{center}
For a preradical $\sigma$, Proposition \ref{preradicales suma directa} implies that $\sigma(\oplus_{i}^{n} \ A_{i})=\oplus_{i=1}^{n} \sigma(A_{i})$ and thus  
\begin{equation*}
\xymatrix{
\oplus_{i=1}^{n} A_{i} \ar[r]^{\ \ \ \ f} & B \\
\oplus_{i=1}^{n} \sigma(A_{i}) \ar[r]^{\ \ \  f|} \ar[u]^{\iota} & \sigma(B)\ar[u]_{\iota}
}
\end{equation*}
is a commutative diagram. Hence, the information defined by the compatible choice assignment $\sigma$ relative to the transformation $f$ is 
\begin{center}
$f(\oplus_{i=1}^{n} \sigma(A_{i}))=f(\sigma(A_{1})) +  f(\sigma(A_{2})) + \dots + f(\sigma(A_{n})) \subseteq \sigma(B)\subseteq B$.
\end{center}

In the above construction, we obtained information on an input-output relationship using a single preradical $\sigma$. We now show that the join operation in $R$-pr allows us to obtain information about the input-output relations regarding a finite set of preradicals. In this case, each component will have an associated preradical with which we will obtain information about the input-output relationship defined by the morphism between the source and the target components.
Thus, the whole construction will gather the amount of information generated by the set of preradicals relative to the functions that comprise the input-output relationships among the different sources and the target component. Suppose we have two morphisms $f_{1}:A_{1}\longrightarrow B$ and $f_{2}:A_{2}\longrightarrow B$ representing an  input-output relationship between the respective sources and target components. 
As we saw above, these morphisms define a morphism from the direct sum $f:A_{1}\oplus A_{2}\longrightarrow B$ such that $f\circ \iota_{i}=f_{i}$ for $i=1,2$. If $\sigma$ and $\tau$ are two preradicals, then we have the commutative diagrams
\begin{equation} \label{diagram 1}
\xymatrix{
A_{1}\ar[r]^{f_{1}} & B \\
\sigma(A_{1}) \ar[r]^{f_{1}|} \ar[u]^{\iota} & \sigma(B)\ar[u]_{\iota}
}
\mbox{   and    } 
\xymatrix{
A_{2}\ar[r]^{f_{2}} & B \\
\tau(A_{2}) \ar[r]^{f_{2}|} \ar[u]^{\iota} & \tau(B)\ar[u]_{\iota}
.}
\tag{*}
\end{equation}
Since $f\circ \iota_{i}=f_{i}$ for $i=1,2$ and $\sigma(B)+\tau(B)\subseteq B$, from (\ref{diagram 1}) we have that 
\begin{equation*}
\xymatrix{
A_{1}\oplus A_{2} \ar[r]^{f} & B \\
\sigma(A_{1})\oplus \tau(A_{2}) \ar[r]^{f|} \ar[u]^{\iota} & \sigma(B)+\tau(B)\ar[u]_{\iota}
}
\end{equation*}
is also a commutative diagram. Note that in the above diagram we are not using directly the preradical $(\sigma \vee \tau)$, although in the right hand we do have $\sigma(B)+\tau(B)=(\sigma \vee \tau)(B)$.

This last argument holds for a number of $n$ inputs coming from components $A_{1},\cdots, A_{n}$ to the component $B$ in a analogous way. In this case, if $\sigma_{1},\cdots \sigma_{n}$ are $n$ preradicals, for each $i\in \{1,\cdots, n\}$ we have the commutative diagram
\begin{equation*}
\xymatrix{
A_{i}\ar[r]^{f_{i}} & B \\
\sigma_{i}(A_{i}) \ar[r]^{f_{i}|} \ar[u]^{\iota} & \sigma_{i}(B)\ar[u]_{\iota}
}
\end{equation*}
which imply the commutativity of the following diagram:
\begin{equation*}
\xymatrix{
\oplus_{i=1}^{n} A_{i} \ar[r]^{\ \ \ \ f} & B \\
\oplus_{i=1}^{n} \sigma_{i}(A_{i}) \ar[r]^{f|} \ar[u]^{\iota} & \Sigma_{i=1}^{n}\sigma_{i}(B)\ar[u]_{\iota}.
}
\end{equation*}
Observe that $\Sigma_{i=1}^{n}\sigma_{i}(B)$ is exactly $(\sigma_{1}\vee \cdots \vee \sigma_{n})(B)$, so as we would expect, the information defined by the set of compatible choice assignments $\sigma_{1},\cdots, \sigma_{n}$ relative to the transformation $f$ is at most $(\sigma_{1}\vee \cdots \vee \sigma_{n})(B)$.

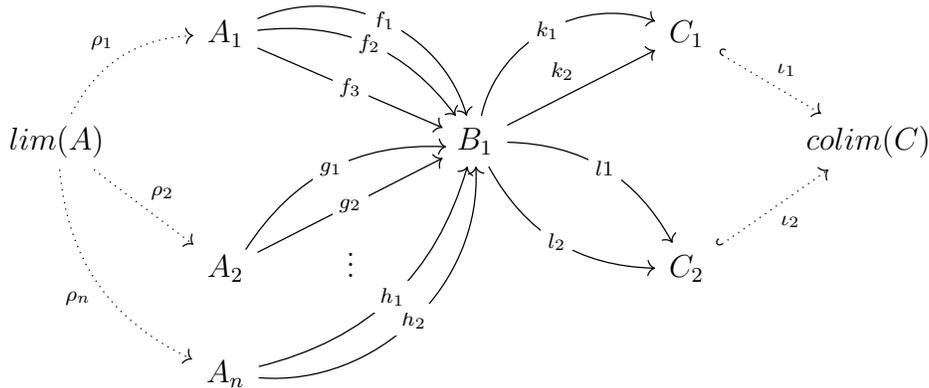
\begin{figure} [h!]
\begin{center}
\begin{tikzcd}
                                                                                                                             & A_{1} \arrow[rrd, "f_{1}" description, bend left=49] \arrow[rrd, "f_{3}" description] \arrow[rrd, "f_{2}" description, bend left] &        &                                                                                                                                                                   &  & C_{1} \arrow[rd, "\iota_{1}", dotted, hook]  &          \\
lim(A) \arrow[ru, "\rho_{1}", dotted, bend left] \arrow[rd, "\rho_{2}", dotted] \arrow[rdd, "\rho_{n}"', dotted, bend right] &                                                                                                                                   &        & B_{1} \arrow[rru, "k_{1}" description, bend left=49] \arrow[rrd, "l{1}" description, bend left] \arrow[rru, "k_{2}"] \arrow[rrd, "l_{2}" description, bend right] &  &                                              & colim(C) \\
                                                                                                                             & A_{2} \arrow[rru, "g_{1}" description, bend left, shift right] \arrow[rru, "g_{2}" description]                                   & \vdots &                                                                                                                                                                   &  & C_{2} \arrow[ru, "\iota_{2}"', dotted, hook] &          \\
                                                                                                                             & A_{n} \arrow[rruu, "h_{2}" description, bend right=49] \arrow[rruu, "h_{1}" description, bend right]                              &        &                                                                                                                                                                   &  &                                              &         
\end{tikzcd}
\caption{Components with multiple connections.}
\label{fig:alpha pre}
\end{center}
\end{figure}

We end this section by illustrating the broader case when the system has multiple distinct input-output relationships between two components. Suppose we have components $A_{1}, \cdots, A_{n}$ whose outputs are the input of a component $B_{1}$, which in turn is also connected to $C_{1}$ and $C_{2}$ (Figure \ref{fig:alpha pre}).

Depending on the information we want to obtain in component $B_{1}$, we may take the limit of $A_{1},\cdots, A_{n}$ which corresponds to their direct product, or we may consider the limit of the whole diagram comprised by $A_{1}, \cdots, A_{n}, B_{1}, C_{1}, C_{2}$ along with their respective input-output relationships (morphisms). The distinction between these two approaches is that in the former one, we gather the plain information coming from the sources $A_{1},\cdots, A_{n}$ whereas in the latter construction, we are considering commutative conditions induced by the limit's property. Since we have already described this last situation in Section \ref{ Persistence}, we will sketch the construction when just the source components define the limit. In this way, $lim(A)=\Pi_{i=1}^{n}\ A_{i}$ and the full information received in component $B_{1}$ is 
\begin{gather*} 
\alpha_{lim(A)}^{lim(A)}(B_{1})=\sum \{\gamma(lim(A)) | \gamma:lim(A)\longrightarrow B_{1} \} \vspace{2mm} \\
=\sum_{i=1}^{3} (f_{i}\circ \rho_{1})(lim(A)) 
+ \sum_{j=1}^{2} (g_{j}\circ \rho_{2})(lim(A)) \vspace{2mm} 
+ \cdots 
+ \sum_{r=1}^{2} (h_{r}\circ \rho_{1})(lim(A)) \vspace{2mm} \\
=\sum_{i=1}^{3} f_{i}(A_{1}) + \sum_{j=1}^{2} g_{j}(A_{2}) + \cdots + \sum_{r=1}^{2} h_{r}(A_{n}).
\end{gather*}
\\ On the other hand, to obtain the information coming from the source components $A_{i}$'s to component $C_{1}$, we take into account 
\begin{center}
$\alpha_{lim(A)}^{lim(A)}(C_{1})=\sum \{\gamma(lim(A)) | \gamma:lim(A)\longrightarrow C_{1} \} \vspace{2mm} $
$=\sum_{i=1}^{3} (k_{1}\circ f_{i}\circ \rho_{1})(lim(A)) + \sum_{j=1}^{2} (k_{1}\circ g_{j}\circ \rho_{2})(lim(A)) + \cdots + \sum_{r=1}^{2} (k_{1}\circ h_{r}\circ \rho_{n})(lim(A)) \vspace{2mm}$
$+ \sum_{i=1}^{3} (k_{2}\circ f_{i}\circ \rho_{1})(lim(A)) + \sum_{j=1}^{2} (k_{2}\circ g_{j}\circ \rho_{2})(lim(A)) + \cdots + \sum_{r=1}^{2} (k_{2}\circ h_{r}\circ \rho_{n})(lim(A))$.
\end{center}
Observe that we can also look at the information in $C_{1}$ that comes from just one of its input-output relationships. For example, we can consider the information coming from the sources that is processed by the morphism $k_{1}$. In this case, we obtain 
\begin{center}
$\sum_{i=1}^{3} (k_{1}\circ f_{i}\circ \rho_{1})(lim(A)) + \sum_{j=1}^{2} (k_{1}\circ g_{j}\circ \rho_{2})(lim(A)) + \cdots + \sum_{r=1}^{2} (k_{1}\circ h_{r}\circ \rho_{n})(lim(A)) \vspace{2mm}$
$=k_{1} \Big( \sum_{i=1}^{3} f_{i}(\rho_{1}(lim(A))) + \sum_{j=1}^{2} g_{j}(\rho_{2}(lim(A))) + \cdots + \sum_{r=1}^{2} h_{r}(\rho_{n}(lim(A))) \Big)$.
\end{center}
This corresponds exactly to the image in $C_{1}$ of the morphism $k_{1}$ restricted to $\alpha_{lim(A)}^{lim(A)}(B_{1})$:
\begin{equation*}
\xymatrix{
B_{1}\ar[r]^{k_{1}} & C_{1} \\
\alpha_{lim(A)}^{lim(A)}(B_{1}) \ar[r]^{k_{1}|} \ar[u]^{\iota} & \alpha_{lim(A)}^{lim(A)}(C_{1})\ar[u]_{\iota}
}
\end{equation*}

The above construction can be used similarly to obtain information of a subset $A_{i_{1}},\cdots, A_{i_{m}}$ of the source components $A_{1},\cdots, A_{n}$ to any other component. For this, we consider the submodule $\Pi_{j=1}^{m} \ A_{i_{j}}$ of $\Pi_{i=1}^{n}\ A_{i}$ to define the preradical $\alpha_{\Pi_{j=1}^{m} \ A_{i_{j}}}^{\Pi_{i=1}^{n}\ A_{i}}$ whose morphisms are those factorized by some $\rho_{i_{l}}$, for $l\in \{1,\dots, m\}$.
\\ Finally, the information that persists from $lim(A)$ through the objects of the system is obtained by considering the $colim(C)$ (i.e., the coproduct of the sink components $C_{1}$ and $C_{2}$) and then evaluating the preradical $\alpha_{lim(A)}^{lim(A)}$ in $colim(C)$: this is the sum of the images of all morphisms $f:lim(A)\longrightarrow colim(C)$ which are induced by the system. This is precisely the full information coming from the source components that flow through the system.

\section{Conclusions} \label{Conclusion}

Preradicals are a tool from category theory that naturally describes the flow of information through a system. Mainly, when we use a $G$-module representation associated with a directed acyclic graph, preradicals characterize the notion of persistence in a $G$-module structure. This description of information's flow can be extended to any quiver representation in $R$-Mod since preradicals are functors that preserve the underlying structure of the representation. Also, by how preradicals are defined, they can be considered as compatible choice assignments that show how certain substructures are related within the modeled system. From a communication theory perspective, this provides a high-order approximation to the information contained in the system. 

\bibliographystyle{amsplain}

\appendix

\section{Homology groups} \label{Apendix}

Given a simplicial complex $X$, a $k$-chain is a formal sum $\Sigma_{i=1}^{n}c_{i} S_{i}$ where each $c_{i}$ is an integer and $S_{i}$ is an oriented $k$-simplex. The collection of $k$-chains forms a group which we denote by $C_{k}(X)$. In fact, $C_{k}(X)$ is a free abelian group whose basis is in one-to-one correspondence with the set of $k$-simplicies of $X$. There is a \textit{boundary} homomorphism
\begin{center}
    $\delta_{k}:C_{k}(X)\longrightarrow C_{k-1}(X)$ 
\end{center}
that  calculates the boundary of a chain and which satisfies $\delta_{k-1}\circ \delta_{k}=0$ for all $k\in \mathbb{N}$. The kernel of $\delta_{k}:C_{k}(X)\longrightarrow C_{k-1}(X)$ is called the $\textit{cycle group}$ and we denote it as $Z_{k}(X)$. The image of 
$\delta_{k+1}:C_{k+1}(X)\longrightarrow C_{k}(X)$ is called the \textit{boundary group} and is denoted by $B_{k}$. Both $Z_{k}$ and $B_{k}$ are contained in $C_{k}(X)$, and since  $\delta_{k-1}\circ \delta_{k}=0$, they satisfy $B_{k} \subseteq Z_{k}(X)$. Thus, the $k$-homology group of $X$ is defined as $H_{k}(X)=Z_{k}/B_{k}$. When the coefficients in the formal sum are taken in a field $F$, then $C_{k}(X),Z_{k}(X),B_{k}(X)$ and $H_{k}(X)$ are all vector spaces. 

Given a filtration $X_{0}\subseteq X_{1} \subseteq \cdots X_{n}$, the persistence homology group $H_{k}^{t}(X_{j})$ is defined as the image of the homomorphism induced by the inclusion map $\iota:X_{j}\hookrightarrow X_{j+t}$. In other words, the persistence homology group $H_{k}^{t}(X_{j})=Im(\iota^{*})$ where $\iota ^{*}:H_{k}(X_{j})\longrightarrow H_{k}(X_{j+t})$. Notice that $H_{k}^{t}(X_{j})$ can be considered as a subgroup of $H_{k}(X_{j+t})$

\end{document}